\def\Date{2008/07/10}


\ifx\pdfoutput\jamaisdefined\else
\input supp-pdf.tex \pdfoutput=1 \pdfcompresslevel=9

\fi

%

\magnification=1200
\hsize=11.25cm
\vsize=18cm
\parskip 0pt
\parindent=12pt
\voffset=1cm
\hoffset=1cm



\catcode'32=9

\font\tenpc=cmcsc10
\font\eightpc=cmcsc8
\font\eightrm=cmr8
\font\eighti=cmmi8
\font\eightsy=cmsy8
\font\eightbf=cmbx8
\font\eighttt=cmtt8
\font\eightit=cmti8
\font\eightsl=cmsl8
\font\sixrm=cmr6
\font\sixi=cmmi6
\font\sixsy=cmsy6
\font\sixbf=cmbx6

\skewchar\eighti='177 \skewchar\sixi='177
\skewchar\eightsy='60 \skewchar\sixsy='60

\catcode`@=11

\def\tenpoint{%
  \textfont0=\tenrm \scriptfont0=\sevenrm \scriptscriptfont0=\fiverm
  \def\rm{\fam\z@\tenrm}%
  \textfont1=\teni \scriptfont1=\seveni \scriptscriptfont1=\fivei
  \def\oldstyle{\fam\@ne\teni}%
  \textfont2=\tensy \scriptfont2=\sevensy \scriptscriptfont2=\fivesy
  \textfont\itfam=\tenit
  \def\it{\fam\itfam\tenit}%
  \textfont\slfam=\tensl
  \def\sl{\fam\slfam\tensl}%
  \textfont\bffam=\tenbf \scriptfont\bffam=\sevenbf
  \scriptscriptfont\bffam=\fivebf
  \def\bf{\fam\bffam\tenbf}%
  \textfont\ttfam=\tentt
  \def\tt{\fam\ttfam\tentt}%
  \abovedisplayskip=12pt plus 3pt minus 9pt
  \abovedisplayshortskip=0pt plus 3pt
  \belowdisplayskip=12pt plus 3pt minus 9pt
  \belowdisplayshortskip=7pt plus 3pt minus 4pt
  \smallskipamount=3pt plus 1pt minus 1pt
  \medskipamount=6pt plus 2pt minus 2pt
  \bigskipamount=12pt plus 4pt minus 4pt
  \normalbaselineskip=12pt
  \setbox\strutbox=\hbox{\vrule height8.5pt depth3.5pt width0pt}%
  \let\bigf@ntpc=\tenrm \let\smallf@ntpc=\sevenrm
  \let\petcap=\tenpc
  \normalbaselines\rm}

\def\eightpoint{%
  \textfont0=\eightrm \scriptfont0=\sixrm \scriptscriptfont0=\fiverm
  \def\rm{\fam\z@\eightrm}%
  \textfont1=\eighti \scriptfont1=\sixi \scriptscriptfont1=\fivei
  \def\oldstyle{\fam\@ne\eighti}%
  \textfont2=\eightsy \scriptfont2=\sixsy \scriptscriptfont2=\fivesy
  \textfont\itfam=\eightit
  \def\it{\fam\itfam\eightit}%
  \textfont\slfam=\eightsl
  \def\sl{\fam\slfam\eightsl}%
  \textfont\bffam=\eightbf \scriptfont\bffam=\sixbf
  \scriptscriptfont\bffam=\fivebf
  \def\bf{\fam\bffam\eightbf}%
  \textfont\ttfam=\eighttt
  \def\tt{\fam\ttfam\eighttt}%
  \abovedisplayskip=9pt plus 2pt minus 6pt
  \abovedisplayshortskip=0pt plus 2pt
  \belowdisplayskip=9pt plus 2pt minus 6pt
  \belowdisplayshortskip=5pt plus 2pt minus 3pt
  \smallskipamount=2pt plus 1pt minus 1pt
  \medskipamount=4pt plus 2pt minus 1pt
  \bigskipamount=9pt plus 3pt minus 3pt
  \normalbaselineskip=9pt
  \setbox\strutbox=\hbox{\vrule height7pt depth2pt width0pt}%
  \let\bigf@ntpc=\eightrm \let\smallf@ntpc=\sixrm
  \let\petcap=\eightpc
  \normalbaselines\rm}
\catcode`@=12

\tenpoint


\long\def\irmaaddress{{%
\bigskip
\eightpoint
\rightline{\quad
\vtop{\halign{\hfil##\hfil\cr
I.R.M.A. UMR 7501\cr
Universit\'e Louis Pasteur et CNRS,\cr
7, rue Ren\'e-Descartes\cr
F-67084 Strasbourg, France\cr
{\tt guoniu@math.u-strasbg.fr}\cr}}\quad}
}}



\catcode`\@=11
\def\pc#1#2|{{\bigf@ntpc #1\penalty \@MM\hskip\z@skip\smallf@ntpc%
	\uppercase{#2}}}
\catcode`\@=12

\def\pointir{\discretionary{.}{}{.\kern.35em---\kern.7em}\nobreak
   \hskip 0em plus .3em minus .4em }

\def\qed{\quad\raise -2pt\hbox{\vrule\vbox to 10pt{\hrule width 4pt
   \vfill\hrule}\vrule}}

\def\rem#1|{\par\medskip{{\it #1}\pointir}}

\def\vspace[#1]{\noalign{\vskip#1}}

\def\abstract#1{\vbox{\eightpoint\narrower\narrower 
\pc ABSTRACT|\pointir #1}}


\def\section#1{\goodbreak\par\vskip .3cm\centerline{\bf #1}
   \par\nobreak\vskip 3pt }

\long\def\th#1|#2\endth{\par\medbreak
   {\petcap #1\pointir}{\it #2}\par\medbreak}

\def\article#1|#2|#3|#4|#5|#6|#7|
    {{\leftskip=7mm\noindent
     \hangindent=2mm\hangafter=1
     \llap{{\tt [#1]}\hskip.35em}{\petcap#2}\pointir
     #3, {\sl #4}, {\bf #5} ({\oldstyle #6}),
     pp.\nobreak\ #7.\par}}
\def\livre#1|#2|#3|#4|
    {{\leftskip=7mm\noindent
    \hangindent=2mm\hangafter=1
    \llap{{\tt [#1]}\hskip.35em}{\petcap#2}\pointir
    {\sl #3}, #4.\par}}
\def\divers#1|#2|#3|
    {{\leftskip=7mm\noindent
    \hangindent=2mm\hangafter=1
     \llap{{\tt [#1]}\hskip.35em}{\petcap#2}\pointir
     #3.\par}}



\catcode`\@=11
\def\c@rr@#1{\vbox{%
  \hrule height \ep@isseur%
   \hbox{\vrule width\ep@isseur\vbox to \t@ille{%
           \vfil\hbox  to \t@ille{\hfil#1\hfil}\vfil}%
            \vrule width\ep@isseur}%
      \hrule height \ep@isseur}}
\def\ytableau#1#2#3#4{\vbox{%
  \gdef\ep@isseur{#2}
   \gdef\t@ille{#1}
    \def\\##1{\c@rr@{$#3 ##1$}}
  \lineskiplimit=-30cm \baselineskip=\t@ille%
    \advance \baselineskip by \ep@isseur%
     \halign{%
      \hfil$##$\hfil&&\kern -\ep@isseur%
       \hfil$##$\hfil \crcr#4\crcr}}}%
\catcode`\@=12

\def\Grille{\setbox13=\vbox to 5mm{\hrule width 110mm\vfill}
\setbox13=\vbox{\offinterlineskip
   \copy13\copy13\copy13\copy13\copy13\copy13\copy13\copy13
   \copy13\copy13\copy13\copy13\box13\hrule width 110mm}
\setbox14=\hbox to 5mm{\vrule height 65mm\hfill}
\setbox14=\hbox{\copy14\copy14\copy14\copy14\copy14\copy14
   \copy14\copy14\copy14\copy14\copy14\copy14\copy14\copy14
   \copy14\copy14\copy14\copy14\copy14\copy14\copy14\copy14\box14}
\ht14=0pt\dp14=0pt\wd14=0pt
\setbox13=\vbox to 0pt{\vss\box13\offinterlineskip\box14}
\wd13=0pt\box13}


\def\fleche(#1,#2)\dir(#3,#4)\long#5{%
\noalign{\nointerlineskip\leftput(#1,#2){\vector(#3,#4){#5}}\nointerlineskip}}


\def\hfl#1#2#3{\smash{\mathop{\hbox to#3{\rightarrowfill}}\limits
^{\scriptstyle#1}_{\scriptstyle#2}}}

\def\gfl#1#2#3{\smash{\mathop{\hbox to#3{\leftarrowfill}}\limits
^{\scriptstyle#1}_{\scriptstyle#2}}}


 \message{`lline' & `vector' macros from LaTeX}
 \catcode`@=11
\def\{{\relax\ifmmode\lbrace\else$\lbrace$\fi}
\def\}{\relax\ifmmode\rbrace\else$\rbrace$\fi}
\def\newcount{\alloc@0\count\countdef\insc@unt}
\def\newdimen{\alloc@1\dimen\dimendef\insc@unt}
\def\newwrite{\alloc@7\write\chardef\sixt@@n}

\newwrite\@unused
\newcount\@tempcnta
\newcount\@tempcntb
\newdimen\@tempdima
\newdimen\@tempdimb
\newbox\@tempboxa

\def\@spaces{\space\space\space\space}
\def\@whilenoop#1{}
\def\@whiledim#1\do #2{\ifdim #1\relax#2\@iwhiledim{#1\relax#2}\fi}
\def\@iwhiledim#1{\ifdim #1\let\@nextwhile=\@iwhiledim
        \else\let\@nextwhile=\@whilenoop\fi\@nextwhile{#1}}
\def\@badlinearg{\@latexerr{Bad \string\line\space or \string\vector
   \space argument}}
\def\@latexerr#1#2{\begingroup
\edef\@tempc{#2}\expandafter\errhelp\expandafter{\@tempc}%
\def\@eha{Your command was ignored.
^^JType \space I <command> <return> \space to replace it
  with another command,^^Jor \space <return> \space to continue without it.}
\def\@ehb{You've lost some text. \space \@ehc}
\def\@ehc{Try typing \space <return>
  \space to proceed.^^JIf that doesn't work, type \space X <return> \space to
  quit.}
\def\@ehd{You're in trouble here.  \space\@ehc}

\typeout{LaTeX error. \space See LaTeX manual for explanation.^^J
 \space\@spaces\@spaces\@spaces Type \space H <return> \space for
 immediate help.}\errmessage{#1}\endgroup}
\def\typeout#1{{\let\protect\string\immediate\write\@unused{#1}}}

\font\tenln    = line10
\font\tenlnw   = linew10

\newdimen\@wholewidth
\newdimen\@halfwidth
\newdimen\unitlength 

\unitlength =1pt


\def\thinlines{\let\@linefnt\tenln \let\@circlefnt\tencirc
  \@wholewidth\fontdimen8\tenln \@halfwidth .5\@wholewidth}
\def\thicklines{\let\@linefnt\tenlnw \let\@circlefnt\tencircw
  \@wholewidth\fontdimen8\tenlnw \@halfwidth .5\@wholewidth}

\def\linethickness#1{\@wholewidth #1\relax \@halfwidth .5\@wholewidth}

\newif\if@negarg

\def\lline(#1,#2)#3{\@xarg #1\relax \@yarg #2\relax
\@linelen=#3\unitlength
\ifnum\@xarg =0 \@vline
  \else \ifnum\@yarg =0 \@hline \else \@sline\fi
\fi}

\def\@sline{\ifnum\@xarg< 0 \@negargtrue \@xarg -\@xarg \@yyarg -\@yarg
  \else \@negargfalse \@yyarg \@yarg \fi
\ifnum \@yyarg >0 \@tempcnta\@yyarg \else \@tempcnta -\@yyarg \fi
\ifnum\@tempcnta>6 \@badlinearg\@tempcnta0 \fi
\setbox\@linechar\hbox{\@linefnt\@getlinechar(\@xarg,\@yyarg)}%
\ifnum \@yarg >0 \let\@upordown\raise \@clnht\z@
   \else\let\@upordown\lower \@clnht \ht\@linechar\fi
\@clnwd=\wd\@linechar
\if@negarg \hskip -\wd\@linechar \def\@tempa{\hskip -2\wd\@linechar}\else
     \let\@tempa\relax \fi
\@whiledim \@clnwd <\@linelen \do
  {\@upordown\@clnht\copy\@linechar
   \@tempa
   \advance\@clnht \ht\@linechar
   \advance\@clnwd \wd\@linechar}%
\advance\@clnht -\ht\@linechar
\advance\@clnwd -\wd\@linechar
\@tempdima\@linelen\advance\@tempdima -\@clnwd
\@tempdimb\@tempdima\advance\@tempdimb -\wd\@linechar
\if@negarg \hskip -\@tempdimb \else \hskip \@tempdimb \fi
\multiply\@tempdima \@m
\@tempcnta \@tempdima \@tempdima \wd\@linechar \divide\@tempcnta \@tempdima
\@tempdima \ht\@linechar \multiply\@tempdima \@tempcnta
\divide\@tempdima \@m
\advance\@clnht \@tempdima
\ifdim \@linelen <\wd\@linechar
   \hskip \wd\@linechar
  \else\@upordown\@clnht\copy\@linechar\fi}

\def\@hline{\ifnum \@xarg <0 \hskip -\@linelen \fi
\vrule height \@halfwidth depth \@halfwidth width \@linelen
\ifnum \@xarg <0 \hskip -\@linelen \fi}

\def\@getlinechar(#1,#2){\@tempcnta#1\relax\multiply\@tempcnta 8
\advance\@tempcnta -9 \ifnum #2>0 \advance\@tempcnta #2\relax\else
\advance\@tempcnta -#2\relax\advance\@tempcnta 64 \fi
\char\@tempcnta}

\def\vector(#1,#2)#3{\@xarg #1\relax \@yarg #2\relax
\@linelen=#3\unitlength
\ifnum\@xarg =0 \@vvector
  \else \ifnum\@yarg =0 \@hvector \else \@svector\fi
\fi}

\def\@hvector{\@hline\hbox to 0pt{\@linefnt
\ifnum \@xarg <0 \@getlarrow(1,0)\hss\else
    \hss\@getrarrow(1,0)\fi}}

\def\@vvector{\ifnum \@yarg <0 \@downvector \else \@upvector \fi}

\def\@svector{\@sline
\@tempcnta\@yarg \ifnum\@tempcnta <0 \@tempcnta=-\@tempcnta\fi
\ifnum\@tempcnta <5
  \hskip -\wd\@linechar
  \@upordown\@clnht \hbox{\@linefnt  \if@negarg
  \@getlarrow(\@xarg,\@yyarg) \else \@getrarrow(\@xarg,\@yyarg) \fi}%
\else\@badlinearg\fi}

\def\@getlarrow(#1,#2){\ifnum #2 =\z@ \@tempcnta='33\else
\@tempcnta=#1\relax\multiply\@tempcnta \sixt@@n \advance\@tempcnta
-9 \@tempcntb=#2\relax\multiply\@tempcntb \tw@
\ifnum \@tempcntb >0 \advance\@tempcnta \@tempcntb\relax
\else\advance\@tempcnta -\@tempcntb\advance\@tempcnta 64
\fi\fi\char\@tempcnta}

\def\@getrarrow(#1,#2){\@tempcntb=#2\relax
\ifnum\@tempcntb < 0 \@tempcntb=-\@tempcntb\relax\fi
\ifcase \@tempcntb\relax \@tempcnta='55 \or
\ifnum #1<3 \@tempcnta=#1\relax\multiply\@tempcnta
24 \advance\@tempcnta -6 \else \ifnum #1=3 \@tempcnta=49
\else\@tempcnta=58 \fi\fi\or
\ifnum #1<3 \@tempcnta=#1\relax\multiply\@tempcnta
24 \advance\@tempcnta -3 \else \@tempcnta=51\fi\or
\@tempcnta=#1\relax\multiply\@tempcnta
\sixt@@n \advance\@tempcnta -\tw@ \else
\@tempcnta=#1\relax\multiply\@tempcnta
\sixt@@n \advance\@tempcnta 7 \fi\ifnum #2<0 \advance\@tempcnta 64 \fi
\char\@tempcnta}

\def\@vline{\ifnum \@yarg <0 \@downline \else \@upline\fi}

\def\@upline{\hbox to \z@{\hskip -\@halfwidth \vrule
  width \@wholewidth height \@linelen depth \z@\hss}}

\def\@downline{\hbox to \z@{\hskip -\@halfwidth \vrule
  width \@wholewidth height \z@ depth \@linelen \hss}}

\def\@upvector{\@upline\setbox\@tempboxa\hbox{\@linefnt\char'66}\raise
     \@linelen \hbox to\z@{\lower \ht\@tempboxa\box\@tempboxa\hss}}

\def\@downvector{\@downline\lower \@linelen
      \hbox to \z@{\@linefnt\char'77\hss}}

\thinlines

\newcount\@xarg
\newcount\@yarg
\newcount\@yyarg
\newcount\@multicnt
\newdimen\@xdim
\newdimen\@ydim
\newbox\@linechar
\newdimen\@linelen
\newdimen\@clnwd
\newdimen\@clnht
\newdimen\@dashdim
\newbox\@dashbox
\newcount\@dashcnt
 \catcode`@=12


\newbox\tbox
\newbox\tboxa

\def\leftzer#1{\setbox\tbox=\hbox to 0pt{#1\hss}%
     \ht\tbox=0pt \dp\tbox=0pt \box\tbox}

\def\rightzer#1{\setbox\tbox=\hbox to 0pt{\hss #1}%
     \ht\tbox=0pt \dp\tbox=0pt \box\tbox}

\def\centerzer#1{\setbox\tbox=\hbox to 0pt{\hss #1\hss}%
     \ht\tbox=0pt \dp\tbox=0pt \box\tbox}

%
\def\image(#1,#2)#3{\vbox to #1{\offinterlineskip
    \vss #3 \vskip #2}}


\def\leftput(#1,#2)#3{\setbox\tboxa=\hbox{%
    \kern #1\unitlength
    \raise #2\unitlength\hbox{\leftzer{#3}}}%
    \ht\tboxa=0pt \wd\tboxa=0pt \dp\tboxa=0pt\box\tboxa}

\def\rightput(#1,#2)#3{\setbox\tboxa=\hbox{%
    \kern #1\unitlength
    \raise #2\unitlength\hbox{\rightzer{#3}}}%
    \ht\tboxa=0pt \wd\tboxa=0pt \dp\tboxa=0pt\box\tboxa}

\def\centerput(#1,#2)#3{\setbox\tboxa=\hbox{%
    \kern #1\unitlength
    \raise #2\unitlength\hbox{\centerzer{#3}}}%
    \ht\tboxa=0pt \wd\tboxa=0pt \dp\tboxa=0pt\box\tboxa}

\unitlength=1mm

\def\cput(#1,#2)#3{\noalign{\nointerlineskip\centerput(#1,#2){#3}
                             \nointerlineskip}}


\ifx\pdfoutput\jamaisdefined
\input epsf.tex

\fi


\parskip 0pt plus 1pt

\def\article#1|#2|#3|#4|#5|#6|#7|
    {{\leftskip=7mm\noindent
     \hangindent=2mm\hangafter=1
     \llap{{\tt [#1]}\hskip.35em}{#2},\quad %
     #3, {\sl #4}, {\bf #5} ({\oldstyle #6}),
     pp.\nobreak\ #7.\par}}
\def\livre#1|#2|#3|#4|
    {{\leftskip=7mm\noindent
    \hangindent=2mm\hangafter=1
    \llap{{\tt [#1]}\hskip.35em}{#2},\quad %
    {\sl #3}, #4.\par}}
\def\divers#1|#2|#3|
    {{\leftskip=7mm\noindent
    \hangindent=2mm\hangafter=1
     \llap{{\tt [#1]}\hskip.35em}{#2},\quad %
     #3.\par}}

\def\l{\lambda}

\def\setT{\mathop{\cal T}}
\def\setB{\mathop{\cal B}}


\rightline{\Date}
\bigskip

\centerline{\bf Hook lengths and shifted parts of partitions}
\bigskip
\centerline{Guo-Niu HAN}
\bigskip\medskip

\abstract{
Some conjectures on partition hook lengths, recently stated by the author,
have been proved and generalized by Stanley, who also needed a formula by
Andrews, Goulden and Jackson on symmetric functions to complete his
derivation. Another identity on symmetric functions can be used instead. The
purpose of this note is to prove it.
}



\def\sec{1}
\section{\sec. Introduction} 
The hook lengths of partitions 
are widely studied in the Theory of
Partitions, in Algebraic Combinatorics and Group Representation Theory.
The basic notions needed here can be found in 
[St99, p.287;  La01, p.1].
A {\it partition}~$\l$ is a sequence of positive 
integers $\l=(\l_1, \l_2,\cdots, \l_\ell)$ such that 
$\l_1\geq \l_2 \geq \cdots \geq \l_\ell>0$.
The integers
$(\l_i)_{i=1,2,\ldots, \ell}$ are called the {\it parts} of~$\l$,
the number $\ell$ of parts being the
{\it length} of $\l$ denoted by $\ell(\l)$.  
The sum of its parts $\l_1+ \l_2+\cdots+ \l_\ell$ is
denoted by $|\l|$.
Let $n$ be an integer, a partition 
$\l$ is said to be a partition of $n$ if $|\l|=n$. We write $\l\vdash n$.
Each partition can be represented by its Ferrers diagram. 
For each box $v$ in the Ferrers diagram of a partition $\l$, or
for each box $v$ in $\l$, for short, define the 
{\it hook length} of $v$, denoted by $h_v(\l)$ or $h_v$, to be the number of 
boxes $u$ such that  $u=v$,
or $u$ lies in the same column as $v$ and above $v$, or in the 
same row as $v$ and to the right of $v$.
The product of all hook lengths of $\l$ is denoted by $H_\l$.

\medskip

The hook length plays an important role in Algebraic Combinatorics 
thanks to the famous hook formula
due to Frame, Robinson and Thrall [FRT54]
$$
f_\l={n!\over H_\l}, \leqno{(\sec.1)}
$$
where $f_\l$ is the number of standard Young tableaux of shape $\l$.
\medskip

For each partition $\l$ let $\l\setminus 1$ be the set of all partitions
$\mu$ obtained from $\l$ by erasing one {\it corner} of $\l$.
By the very construction of the standard Young tableaux and (\sec.1) we have
$$
f_\l=\sum_{\mu\in \l\setminus 1} f_{\mu}  \leqno{(\sec.2)}
$$
and then
$$
{n\over H_\l}=\sum_{\mu\in \l\setminus 1} {1 \over H_{\mu}}.  \leqno{(\sec.3)}
$$

\medskip
In this note we establish the following perturbation of formula (\sec.3).
Define the {\it $g$-function} of a partition $\l$ of $n$ to be
$$
g_\l(x)= \prod_{i=1}^n(x+\l_i-i),  \leqno{(\sec.4)}
$$
where $\l_i=0$ for $i\geq \ell(\l)+1$.

\proclaim Theorem \sec.1. 
Let $x$ be a formal parameter. For each partition $\l$ we have
$$
{ g_\l(x+1) - g_\l(x) \over H_\l}
=
\sum_{\mu\in\l\setminus 1} {g_\mu (x) \over H_\mu }.
\leqno{(\sec.5)}
$$

Theorem \sec.1 is proved in Section 2.
Some equivalent forms of Theorem \sec.1 and remarks are given in Section 4.
As an application we prove (see Section 3) 
the following result due to Stanley [St08].

\proclaim Theorem \sec.2.
Let $p, e$ and $s$ be the usual symmetric functions [Ma95, Chap.I]. Then 
$$
\sum_{k=0}^n {x+k-1\choose k} p_1^k e_{n-k}
=\sum_{\l\vdash n} H_\l^{-1} g_\l(x+n) s_\l. 
\leqno{(\sec.6)}
$$

Recently, the author stated some conjectures on partition hook lengths [Ha08a],
which were suggested by hook length expansion techniques (see [Ha08b]).
Later, Conjecture 3.1 in [Ha08a] was proved by Stanley  [St08]. One step of
his proof is formula (\sec.6), based on 
a result by Andrews, Goulden and Jackson [AGJ88]. 
In this paper we provide a simple and direct proof of formula (\sec.6). 
\medskip

{\it Remark}. Let $D$ be the difference operator defined by 
$$D(f(x)) = f(x+1)-f(x).$$
By iterating formula (\sec.5) we obtain
$$
D^{n} {g_\l(x)\over H_\l} = f_\l,
$$
which is precisely the hook length formula (\sec.1).

\def\sec{2}
\section{\sec. Proof of Theorem 1.1} 

Let 
$$
\epsilon(x)=
{ g_\l(x+1) - g_\l(x) \over H_\l}
-
\sum_{\mu\in\l\setminus 1} {g_\mu (x) \over H_\mu }.
\leqno{(\sec.1)}
$$
We see that $\epsilon(x)$ is a polynomial in $x$ whose degree is 
less than or equal to $n$.  Moreover 
$$
[x^n] \epsilon(x) = 
[x^n] { g_\l(x+1) - g_\l(x) \over H_\l} =0.
$$
Furthermore, 
$$
[x^{n-1}] g_\l(x+1)  = \sum_{i=1}^n (\l_i-i+1) = n+ \sum_{i=1}^n (\l_i-i)
=n+[x^{n-1}] g_\l(x)  
$$
and 
$$
[x^{n-1}] \epsilon(x) = 
[x^{n-1}] { g_\l(x+1) - g_\l(x) \over H_\l}  -
\sum_{\mu\in\l\setminus 1} {1 \over H_\mu }
=
 { n \over H_\l}  -
\sum_{\mu\in\l\setminus 1} {1 \over H_\mu }
=0.
$$
The last equality is guaranteed by (1.3),
so that $\epsilon(x)$ is a polynomial in $x$ whose degree is 
less than and equal to $n-2$.  To prove that $\epsilon(x)$ is actually zero, 
it suffices 
to find $n-1$ distinct values for $x$
such that $\epsilon(x)=0$. 
In the following we prove that $\epsilon(i-\l_i)=0$ for 
$i-\l_i$ for $i=1,2,\ldots, n-1$. 

\medskip

If $\l_i = \l_{i+1}$, or if
the $i$-th row has no corner, the factor $x+\l_i-i$ lies 
in $g_\l(x)$ and also in $g_\mu(x)$ for all $\mu\in\l\setminus 1$ . 
The factor 
$(x+1)+\l_{i+1} - (i+1) =x+\l_{i} - i$ is furthermore in $g_\l(x+1)$, 
so that 
$\epsilon(i-\l_i)=0$.
\medskip

Next, if $\l_i \geq \l_{i+1}+1$, or if the $i$-th row has a corner,
the factor $x+\l_i-i$
lies in $g_\l(x)$ and $g_\mu(x)$ for all $\mu\in\l\setminus 1$,
except for $\mu=\l'$, which is the partition obtained from $\l$ by erasing the
corner from the $i$-th row. In this case equality (\sec.1) becomes
$$
\epsilon(i-\l_i)={g_\l(i-\l_i+1) \over H_\l} - {g_{\l'}(i-\l_i)\over H_{\l'}}.
$$
For proving Theorem 1.1, it remains to prove $\epsilon(i-\l_i)=0$ or
$$
{H_\l\over H_{\l'}}={g_\l(i-\l_i+1) \over g_{\l'}(i-\l_i) }. \leqno{(\sec.2)} 
$$

Consider the following product
$$
{g_\l(x+1) \over g_{\l'}(x)} = 
{\prod_{j=1}^n(x+\l_j-j+1) \over \prod_{j=1}^{n-1}(x+\l'_j-j) }. 
\leqno{(\sec.3)}
$$

The set of all $1\leq j \leq n-1$ such that $\l_j>\l_{j+1}$ is denoted by
$\setT$. For $1\leq j\leq n-1$ and  $j\not\in\setT$ 
(which implies that $j\not=i$ and $\l'_j=\l_j=\l_{j+1}$), the numerator
contains $x+\l_{j+1}-(j+1)+1 = x+\l_j-j$ and the denominator also contains
$x+\l'_j-j=x+\l_j-j$. After cancellation of those common factors, 
(\sec.3) becomes

$$
{g_\l(x+1) \over g_{\l'}(x)} = 
{\prod_{j\in\setB}(x+\l_j-j+1) 
	\over 
\prod_{j\in\setT}(x+\l'_j-j) } 
\leqno{(\sec.4)}
$$
where $\setB=\{1\} \cup \{i+1 \mid i\in\setT\}$.
Letting $x=i-\l_i$ in (\sec.4) yields

$$
{g_\l(i-\l_i+1) \over g_{\l'}(i-\l_i)} = 
{\prod_{j\in\setB}(i-\l_i+\l_j-j+1) 
	\over 
\prod_{j\in\setT}(i-\l_i+\l'_j-j) }.
\leqno{(\sec.5)}
$$

We distinguish the factors in the right-hand side of (\sec.5) as follows.

(C1) For $j\in\setB$ and $j>i$,
$
i-\l_i+\l_j-j+1 = 
-(\l_i-\l_j+j-i-1) = -h_v(\l),
$
where $v$ is the box $(i, \l_j+1)$ in $\l$.

(C2) For $j\in\setB$ and $j\leq i$,
$
i-\l_i+\l_j-j+1 =  h_v(\l),
$
where $v$ is the box $(j, \l_i)$ in $\l$.

(C3) For $j\in\setT$ and $j>i$,
$
i-\l_i+\l_j-j = 
-(\l_i-\l_j+j-i) = -h_u(\l'),
$
where $u$ is the box $(i, \l_j)$ in $\l'$.

(C4) For $j\in\setT$ and $j< i$,
$
i-\l_i+\l_j-j =  h_u(\l'),
$
where $u$ is the box $(j, \l_i)$ in $\l'$.

(C5) For $j\in\setT$ and $j= i$,
$
i-\l_i+\l_j'-j =  
i-\l_i+\l_i'-i = -1. 
$
See Fig. \sec.3 and \sec.4 for an example.

Since each $j\in \setB$ such that $j>i$ is associated with $j-1\in\setT$ and $j-1\geq i$, the
right-hand side of (\sec.5) is positive and can be re-written

$$
{g_\l(i-\l_i+1) \over g_{\l'}(i-\l_i)} = 
{\prod_{v}h_v(\l)
	\over 
\prod_{u}h_u(\l') } ,
\leqno{(\sec.6)}
$$
where $v,u$ range over the boxes described in (C1)-(C4).
Finally $H_\l / H_{\l'}$ is equal to the right-hand side of
(\sec.6), since the hook lengths of all other boxes cancel.
We have completed the proof of (\sec.2). \qed

\medskip

For example, consider the partition $\l=55331$ and $i=4$. We have 
$\l'=55321$ and
$$
{H_\l \over H_{\l'}} = 
{4\cdot 2\cdot 1\cdot 2\cdot 5\cdot 6 \over 3\cdot 1\cdot 1\cdot 4\cdot 5}.
=
{4\cdot 2\cdot 2\cdot  6 \over 3\cdot 4}.
$$

{
\long\def\maplebegin#1\mapleend{}

\maplebegin

# --------------- begin maple ----------------------

# Copy the following text  to "makefig.mpl"
# then in maple > read("makefig.mpl");
# it will create a file "z_fig_by_maple.tex"

#\unitlength=1pt

Hu:= 12.4; # height quantities
Lu:= Hu; # large unity

X0:=-95.0; Y0:=15.6; # origin position

File:=fopen("z_fig_by_maple.tex", WRITE);

mhook:=proc(x,y,lenx, leny)
local i, d,sp, yy, xx, ct;
	sp:=Hu/8;
	ct:=0;
	for xx from x*Lu+X0 to x*Lu+X0+Hu by sp do
		yy := y*Hu+Y0; 
		fprintf(File, "\\vline(
				xx,   yy+sp*ct-0.2,    Lu*leny-sp*ct+0.1);
		ct:=ct+1;
	od:
	
	ct:=0;
	for yy from y*Hu+Y0 to y*Hu+Y0+Hu by sp do
		xx := x*Lu+X0; 
		fprintf(File, "\\hline(
				xx+sp*ct-0.2,   yy,    Lu*lenx-sp*ct+0.1);
		ct:=ct+1;
	od:

end;

mhook(1,1,2,3);

fclose(File);
# -------------------- end maple -------------------------
\mapleend

\setbox1=\hbox{$
\def\b{\\{\hbox{}}}
\ytableau{12pt}{0.4pt}{}
{
 \b           \cr
 \\4 &\\2 &\\1   \cr 
 \b &\b &\\2   \cr
 \b &\b &\\5 &\b &\b   \cr
 \b &\b &\\6 &\b &\b   \cr
\noalign{\vskip 3pt}
\noalign{\hbox{Fig.~\sec.1. Hook lengths of $\l$}}
}$
}
\setbox2=\hbox{$
\def\b{\\{\hbox{}}}
\def\oc#1{\,#1}
\ytableau{12pt}{0.4pt}{}
{
 \b           \cr
 \\3 &\\1    \cr 
 \b &\b &\\1   \cr
 \b &\b &\\4 & \b & \b     \cr
 \b &\b &\\5 &\b &\b   \cr
\noalign{\vskip 3pt}
\noalign{\hbox{Fig.~\sec.2. Hook lengths of $\l'$}}
}$
}
$$\box1\qquad\qquad\box2$$
}

On the other hand, $\setT=\{2,4,5\}$, $\setB=\{1,3,5,6\}$ and
$$
{g_\l(x+1) \over g_{\l'}(x)} = 
{ (x+5) (x+1) (x-3) (x-5) \over (x+3) (x-2) (x-4) }.
$$
Letting $x=i-\l_i=4-3=1$ yields
$$
{g_\l(2) \over g_{\l'}(1)} = 
{ (6) (2) (-2) (-4) \over (4) (-1) (-3) }
= 
{ 6\cdot 2\cdot 2\cdot 4 \over 4\cdot 3 }.
$$

{
\long\def\maplebegin#1\mapleend{}

\maplebegin

# --------------- begin maple ----------------------

# Copy the following text  to "makefig.mpl"
# then in maple > read("makefig.mpl");
# it will create a file "z_fig_by_maple.tex"

#\unitlength=1pt

Hu:= 12.4; # height quantities
Lu:= Hu; # large unity

X0:=-95.0; Y0:=15.6; # origin position

File:=fopen("z_fig_by_maple.tex", WRITE);

mhook:=proc(x,y,lenx, leny)
local i, d,sp, yy, xx, ct;
	sp:=Hu/8;
	ct:=0;
	for xx from x*Lu+X0 to x*Lu+X0+Hu by sp do
		yy := y*Hu+Y0; 
		fprintf(File, "\\vline(
				xx,   yy+sp*ct-0.2,    Lu*leny-sp*ct+0.1);
		ct:=ct+1;
	od:
	
	ct:=0;
	for yy from y*Hu+Y0 to y*Hu+Y0+Hu by sp do
		xx := x*Lu+X0; 
		fprintf(File, "\\hline(
				xx+sp*ct-0.2,   yy,    Lu*lenx-sp*ct+0.1);
		ct:=ct+1;
	od:

end;

mhook(1,1,2,3);

fclose(File);
# -------------------- end maple -------------------------
\mapleend

\setbox1=\hbox{$
\def\b{\\{\hbox{}}}
\ytableau{12pt}{0.4pt}{}
{
 \b           \cr
 \\v &\\v &\\1   \cr 
 \b &\b &\\v   \cr
 \b &\b &\\X &\b &\b   \cr
 \b &\b &\\v &\b &\b   \cr
\noalign{\vskip 3pt}
\noalign{\hbox{Fig.~\sec.3. The boxes $v$ in $\l$}}
}$
}
\setbox2=\hbox{$
\def\b{\\{\hbox{}}}
\def\oc#1{\,#1}
\ytableau{12pt}{0.4pt}{}
{
 \b           \cr
 \\u &\\1    \cr 
 \b &\b &\\1   \cr
 \b &\b &\\u & \b & \b     \cr
 \b &\b &\\X &\b &\b   \cr
\noalign{\vskip 3pt}
\noalign{\hbox{Fig.~\sec.4. The boxes $u$ in $\l'$}}
}$
}
$$\box1\qquad\qquad\box2$$
}


\def\sec{3}
\section{\sec. Proof of Theorem 1.2} 

Let $R_n(x)$ be the right-hand side of (1.6). By Theorem 1.1
$$
\leqalignno{
R_n(x) &=\sum_{\l\vdash n} \Bigl(
	{ g_\l(x+n-1) \over H_\l} + \sum_{\mu\in\l\setminus 1}
	{ g_\mu(x+n-1) \over H_\mu} \Bigr) s_\l \cr
&= R_n(x-1)+ 
	\sum_{\l\vdash n} 
	\sum_{\mu\in\l\setminus 1}
	{ g_\mu(x+n-1) \over H_\mu}  s_\l \cr
&= R_n(x-1)+ 
	\sum_{\mu\vdash n-1} 
	\sum_{\l\,:\, \mu\in\l\setminus 1}
	{ g_\mu(x+n-1) \over H_\mu}  s_\l\cr
&= R_n(x-1)+ 
	\sum_{\mu\vdash n-1} 
	{ g_\mu(x+n-1) \over H_\mu}  p_1 s_\mu,\cr
}
$$
where the next to last equality is 
$$
	\sum_{\l\, :\, \mu\in\l\setminus 1} s_\l = p_1 s_\mu
$$
by using Pieri's rule [Ma95, p.73].
We obtain the following recurrence for $R_n(x)$.
$$
R_n(x)=R_{n}(x-1) + p_1 R_{n-1}(x). \leqno{(\sec.1)}
$$

Let $L_n(x)$ be the left-hand side of (1.6). Using elementary 
properties of binomial coefficients
$$
\leqalignno{
L_n(x)&=\sum_{k=0}^n {x+k-1\choose k} p_1^k e_{n-k} \cr
&=  e_{n} +  \sum_{k=1}^n ({x+k-2\choose k}+{x+k-2\choose k-1}) 
   p_1^k e_{n-k} \cr
&=  L_n(x-1)  +  p_1\sum_{k=1}^n {x+k-2\choose k-1} p_1^{k-1} e_{n-k} \cr
&=  L_n(x-1)  +  p_1 L_{n-1} (x). &(\sec.2) \cr
}
$$
We verify that $L_1(x)=R_1(x)$ and $L_n(0)=R_n(0)$, so that
$L_n(x)=R_n(x)$ by (\sec.1) and (\sec.2). \qed

\def\sec{4}
\section{\sec. Equivalent forms and further remarks} 

Let $\l=\l_1\l_2\cdots\l_\ell$ be a partition of $n$.
The set of all $1\leq j \leq n$ such that $\l_j>\l_{j+1}$ is denoted by
$\setT$ and let $\setB=\{1\} \cup \{i+1 \mid i\in\setT\}$.
Those two sets can be viewed as the {\it in-corner} and {\it out-corner} 
index sets, respectively. Notice that
$\#\setB = \#\setT +1$.
For each $i\in\setT$ we define $\l^{i-}$ to be the partition of $n-1$
obtained form $\l$ by erasing the right-most box from the $i$-th row. Hence 
$$
\l\setminus 1 =\{\l^{i-} \mid i\in \setT\}. \leqno{(\sec.1)}
$$
We verify that
$$
g_{\l^{i-}}(x)={ g_\l(x) (x+\l_i-i-1)\over (x+\l_i -i)(x-n)}.
\leqno{(\sec.2)}
$$
From Theorem 1.1
$$
{ g_\l(x+1) - g_\l(x) \over H_\l}
=
\sum_{i\in\setT} { g_\l(x) (x+\l_i-i-1)\over (x+\l_i -i)(x-n)}
{1 \over H_\l^{i-} }
$$
or
$$
\sum_{i\in\setT} {H_\l\over H_{\l^{i-}}} \times (1-{1\over x+\l_i -i})=
		n-x+{(x-n) g_\l(x+1) \over g_\l(x)}.
		\leqno{(\sec.3)}
$$
Let us re-write (1.3) 
$$
\sum_{\mu\in\l\setminus 1} {H_\l\over H_\mu} =n. \leqno{(\sec.4)} 
$$
By subtracting (\sec.3) from (\sec.4) we obtain the following
equivalent form of Theorem 1.1.

\proclaim Theorem \sec.1. 
We have
$$
\sum_{i\in\setT} {H_\l\over H_{\l^{i-}}} \times{1\over x+\l_i -i}=
		x-{(x-n) g_\l(x+1) \over g_\l(x)}. 
\leqno{(\sec.5)}
$$

By the definitions of $\setT$ and $\setB$ we have
$$
{(x-n) g_\l(x+1) \over g_\l(x)}
={\prod_{i\in\setB }(x+\l_i-i+1) \over \prod_{i\in\setT}(x+\l_i-i) },
\leqno{(\sec.6)}
$$
so that Theorem 1.1 is also equivalent to the following result.

\proclaim Theorem \sec.2. 
We have
$$
\sum_{i\in\setT} {H_\l\over H_{\l^{i-}}} \times{1\over x+\l_i -i}=
		x-{\prod_{i\in\setB }(x+\l_i-i+1) \over \prod_{i\in\setT}(x+\l_i-i) }. 
\leqno{(\sec.7)}
$$

For example, take $\l=55331$. Then $\setT={2,4,5}$ and $\setB={1,3,5,6}
=\{1, 2+1,4+1,5+1\}$.
{
\long\def\maplebegin#1\mapleend{}

\maplebegin

# --------------- begin maple ----------------------

# Copy the following text  to "makefig.mpl"
# then in maple > read("makefig.mpl");
# it will create a file "z_fig_by_maple.tex"

#\unitlength=1pt

Hu:= 12.4; # height quantities
Lu:= Hu; # large unity

X0:=-95.0; Y0:=15.6; # origin position

File:=fopen("z_fig_by_maple.tex", WRITE);

mhook:=proc(x,y,lenx, leny)
local i, d,sp, yy, xx, ct;
	sp:=Hu/8;
	ct:=0;
	for xx from x*Lu+X0 to x*Lu+X0+Hu by sp do
		yy := y*Hu+Y0; 
		fprintf(File, "\\vline(
				xx,   yy+sp*ct-0.2,    Lu*leny-sp*ct+0.1);
		ct:=ct+1;
	od:
	
	ct:=0;
	for yy from y*Hu+Y0 to y*Hu+Y0+Hu by sp do
		xx := x*Lu+X0; 
		fprintf(File, "\\hline(
				xx+sp*ct-0.2,   yy,    Lu*lenx-sp*ct+0.1);
		ct:=ct+1;
	od:

end;

mhook(1,1,2,3);

fclose(File);
# -------------------- end maple -------------------------
\mapleend

\setbox1=\hbox{$
\def\b{\\{\hbox{}}}
\ytableau{12pt}{0.4pt}{}
{
 \\A           \cr
 \b &\b &\\B   \cr 
 \b &\b &\b   \cr
 \b &\b &\b &\b &\\C   \cr
 \b &\b &\b &\b &\b   \cr
\noalign{\vskip 3pt}
\noalign{\hbox{Fig.~\sec.1. in-corner}}
}$
}
\setbox2=\hbox{$
\def\b{\\{\hbox{}}}
\def\oc#1{\,#1}
\ytableau{12pt}{0.4pt}{}
{
 \oc{a}            \cr
 \b &\oc{b}          \cr
 \b &\b &\b   \cr 
 \b &\b &\b &\oc{c}  \cr
 \b &\b &\b &\b &\b   \cr
 \b &\b &\b &\b &\b &\oc{d}  \cr
\noalign{\vskip 3pt}
\noalign{\hbox{Fig.~\sec.2. out-corner}}
}$
}
$$\box1\qquad\qquad\box2$$
}
Hence
$\l^{2-}=54331$, $\l^{4-}=55321$ and $\l^{5-}=55330$. 
Equality (\sec.7) becomes
$$
\leqalignno{
&{H_\l \over {H_{\l^{5-}}}} \times {1\over x-4}
+{H_\l \over {H_{\l^{4-}}}} \times {1\over x-1}
+{H_\l \over {H_{\l^{2-}}}} \times {1\over x+3} \cr
&\qquad\qquad\qquad\qquad=x-{(x-5)(x-3)(x+1)(x+5)\over (x-4)(x-1)(x+3) }.\cr
&\qquad\qquad\qquad\qquad={17x^2 -38x-75 \over (x-4)(x-1)(x+3) }.\cr
}
$$

Theorems \sec.1 and \sec.2 can be proved directly using the method used in
the proof of Theorem 1.1. First, we must verify that the numerator in
the right-hand side of (\sec.5) is a polynomial in $x$ whose
degree is less than ($\leq$) $\#\setT-1$. By the partial fraction expansion
technique it suffices to verify that (\sec.7) is true for all 
$x=i-\l_i$ ($i\in\setT$). This direct proof contains the main part
of the proof of Theorem 1.1. However it does not make use of the fundamental 
relation (1.3) or (\sec.4). Thus, the following corollary of Theorem \sec.2
makes sense.

\proclaim Corollary \sec.4.
We have
$$
\sum_{\mu\in\l\setminus 1} {H_\l\over H_\mu} =n. \leqno{(\sec.8)} 
$$

{\it Proof}.
Let $\#\setT=k$. The right-hand side of (\sec.7)
has the following form	
$$
{Cx^{k-1} + \cdots \over x^k + \cdots}.
$$
We now evaluate the coefficient $C$. By (4.6) we can write  $C=A-B$ with
$$
A=[x^{n-1}] x\prod_{i=1}^n (x+\l_i-i) = \sum_{1\leq i<j\leq n} (\l_i-i)(\l_j-j)
$$
and 
$$
\leqalignno{
B&=[x^{n-1}] (x-n)\prod_{i=1}^n (x+\l_i-i+1) \cr
&= \sum_{1\leq i<j\leq n} (\l_i-i+1)(\l_j-j+1) 
-n \sum_{1\leq i\leq n} (\l_i-i+1).\cr
&= B_1 -n \sum_{1\leq i\leq n} (\l_i-i+1),\cr
}
$$
where
$$
\leqalignno{
B_1&= \sum_{1\leq i<j\leq n} (\l_i-i+1)(\l_j-j+1) \cr
&= \sum_{1\leq i<j\leq n} 
\Bigl((\l_i-i)(\l_j-j) + (\l_i-i) +(\l_j-j) +1 \Bigr)  \cr
&= A+ \sum_{1\leq i<j\leq n}  (\l_i-i) 
+\sum_{1\leq i<j\leq n} (\l_j-j) +{n\choose 2}  \cr
&= A+ \sum_{1\leq i\leq n}  (n-i)(\l_i-i) 
+\sum_{1\leq j\leq n} (j-1)(\l_j-j) +{n\choose 2}  \cr
&= A+ \sum_{1\leq i\leq n}  (n-1)(\l_i-i) +{n\choose 2}.  \cr
}
$$
Finally
$$
\leqalignno{
C&=A-B\cr
	&=
- \sum_{1\leq i\leq n}  (n-1)(\l_i-i) -{n\choose 2}
+n \sum_{1\leq i\leq n} (\l_i-i+1)\cr
&=
- \sum_{1\leq i\leq n}  n(\l_i-i) 
+ \sum_{1\leq i\leq n}  (\l_i-i) 
-{n\choose 2}
+n \sum_{1\leq i\leq n} (\l_i-i)
+n^2 \cr
&=
\sum_{1\leq i\leq n}  (\l_i-i) 
-{n\choose 2}
+n^2 \cr
&=
n- {n+1\choose 2}
-{n\choose 2}
+n^2 = n. \qed \cr
}
$$


\bigskip \bigskip


\centerline{References}

{\eightpoint

\bigskip 
\bigskip

\article AGJ88|Andrews, G; Goulden, I; Jackson, D. M.|Generalizations of Cauchy's summation formula for Schur functions|Trans. Amer. Math. Soc.|310|1988|805--820|

\article FRT54|Frame, J. Sutherland; Robinson, Gilbert de Beauregard;        
Thrall, Robert M.|The hook graphs of the symmetric groups|Canadian 
J. Math.|6|1954|316--324|

\divers Ha08a|Han, Guo-Niu|\quad Some conjectures and open problems
about partition hook length, {\sl Experimental Mathematics}, in press,
15 pages, {\oldstyle 2008}|

\divers Ha08b|Han, Guo-Niu|Discovering hook length formulas by 
expansion technique, {\sl in preparation}, 42 pages, {\oldstyle 2008}|

\livre La01|Lascoux, Alain|Symmetric Functions and Combinatorial Operators on 
Polynomials|CBMS Regional Conference Series in Mathematics, Number 99, 
{\oldstyle 2001}|

\livre Ma95|Macdonald, Ian G.|Symmetric Functions and Hall Polynomials|
Second Edition, Clarendon Press, Oxford, {\oldstyle 1995}|

\livre St99|Stanley, Richard P.|Enumerative Combinatorics, vol. 2|
Cambridge University Press, {\oldstyle 1999}|

\divers St08|Stanley, Richard P.|Some combinatorial properties of hook 
lengths, contents, and parts of partitions, 
{\bf arXiv:0807.0383 [math.CO]},
18 pages	, {\oldstyle 2008}|

\bigskip

\irmaaddress
}
\vfill\eject

\end